\documentclass{article}
\usepackage{amssymb}

\def\dbigcup{\bigcup}
\begin{document}

\title{Approximation on the complex sphere }
\author{Huda Alsaud, Alexander Kushpel\thanks{%
This research has been supported by the EPSRC Grant EP/H020071/1.}, Jeremy
Levesley \\
Department of Mathematics\\
University of Leicester, UK\\
E-mail: huda\_leen@hotmail.com, ak412@le.ac.uk, jl1@le.ac.uk}
\date{25 June 2011}
\maketitle

\begin{abstract}
We develop new elements of harmonic analysis on the complex sphere on the
basis of which Bernstein's, Jackson's and Kolmogorov's inequalities are
established. We apply these results to get order sharp estimates of $m$-term
approximations. The results obtained is a synthesis of new results on
classical orthogonal polynomials, harmonic analysis on manifolds and
geometric properties of Euclidean spaces.
\end{abstract}

\bigskip

Keywords: approximation, complex sphere, volume, $m$-term approximation 

Subject: 46B06, 42B15, 43A90, 33C50 



\section{Introduction}

Let $X$ be a separable real Banach space and $\Xi :=\{\xi _{k}\}_{k\in 
\mathbb{N}}$ be a dense subset of $X$, i.e., $\mathsf{cl}_{X}\left( \Xi
\right) =X$. For a fixed $m\in \mathbb{N}$ let $\Omega _{m}:=\{k_{1}<\cdots
<k_{m}\}\subset \mathbb{N}$ and $\Xi (\Omega _{m}):=\mathrm{lin}\{\xi
_{k_{l}}\}_{l=1}^{m}$. Consider the best approximation of an element $\phi
\in X$ by the subspace $\Xi (\Omega _{m})$ in $X$, 
\[
\nu (\phi ,\Xi (\Omega _{m}),X):=\inf_{\xi \in \Xi (\Omega _{m})}\,\,\Vert
\phi -\xi \Vert _{X}=\inf_{(\alpha _{1},\cdots ,\alpha _{m})\in \mathbb{R}%
^{m}}\,\,\left\Vert \phi -\sum_{l=1}^{m}\alpha _{l}\,\xi _{k_{l}}\right\Vert
_{X} 
\]%
The best $m$-term approximation of $\phi \in X$ with regard to the given
system $\Xi $ (frequently $\Xi $ is called dictionary) is 
\[
\nu _{m}(\phi ,\Xi ,X):=\inf_{\Omega _{m}\subset \mathbb{N}}\,\,\,\nu (\phi
,\Xi (\Omega _{m}),X). 
\]%
Finally, $m$-term approximation of a given set $\mathcal{K}\subset X$ is 
\[
\nu _{m}:=\nu _{m}(\mathcal{K},\Xi ,X):=\sup_{\phi \in \mathcal{K}}\,\,\,\nu
_{m}(\phi ,\Xi ,X). 
\]%
$m$-Term approximation has been introduced by Stechkin \cite{stechkin} in
the case $X=L_{p}\left( \mathbb{S}^{1}\right) ,$ $p=2$, then studied by
Ismagilov \cite{ismagilov} and many others for any $1\leq p\leq \infty .$
Here, $\mathbb{S}^{1}$ is the unit circle. During the last years $m$-term
approximations and $n$-widths became very popular in numerical methods for
PDE's. More specifically, in recently developed reduced basis methods. Also,
the idea of so-called "greedy algorithms" has been inspired by $m$-term
approximations. It is natural to call $m$-term approximations considered
here as harmonic $m$-widths by analogy with known trigonometric $m$-widths.
Remark that Kolmogorov's $n$-widths, defined as 
\[
d_{n}\left( \mathcal{K},X\right) :=\inf_{L_{n}\subset X}\,\,\sup_{x\in 
\mathcal{K}}\,\,\inf_{y\in L_{n}}\,\Vert x-y\Vert _{X}, 
\]%
where $\mathcal{K}$ is a centrally symmetric compact in $X$, can be bigger,
less or equal to the respective $n$-term approximations. Observe that $m$%
-term approximation is a highly nonlinear method of approximation. In
particular, in this article we show that in the case of Sobolev's classes $%
W_{p}^{\gamma }$ it is not possible to improve the rate of convergence in $%
L_{q}$, $1\leq q\leq p\leq \infty $ using $m$-term approximation instead of
linear polynomial approximation.

Our lower bounds of $m-$term approximations are essentially based on
Bernstein's inequality \cite{lk} 
\begin{equation}
\Vert t_{N}^{(\gamma )}\Vert _{q}\leq N^{\gamma +d(1/p-1/q)_{+}}\Vert
t_{N}\Vert _{p},\gamma >0,\,\,1\leq p,q\leq \infty ,\forall t_{N}\in 
\mathcal{T}_{N},  \label{ber}
\end{equation}%
where $\mathcal{T}_{N}$ is defined in (\ref{1111}) and methods of Geometry
of Banach spaces. We will need some general definitions.

Let $\Omega _{m}:=\{k_{1}<\cdots <k_{m}\}\subset \mathbb{N}$ and $\Xi
_{n}(\Omega _{m}):=\mathrm{lin}\{\mathrm{H}_{k_{l}}\}_{l=1}^{m},$ where $%
n=\dim \mathrm{lin}\{\mathrm{H}_{k_{l}}\}_{l=1}^{m}$ and $\mathrm{H}_{k_{l}}$
is an eigenspace of Laplace-Beltrami operator on the complex sphere $\mathbb{%
S}^{d}(\mathbb{C})$ defined in the Section \ref{harmonic analysis}. In the
special case $\Omega _{N}=\left\{ 1,2,\cdot \cdot \cdot ,N\right\} $ we
shall write 
\begin{equation}  \label{1111}
\mathcal{T}_{N}:=\mathrm{lin}\{\mathrm{H}_{k}\}_{k=1}^{N}.
\end{equation}

Let $\{\xi _{k}\}_{k\in \mathbb{N}}$ be a sequence of orthonormal, functions
on $\mathbb{S}^{d}(\mathbb{C})$. Let $X$ be a Banach space of functions on $%
\mathbb{S}^{d}(\mathbb{C})$ with the norm $\Vert \cdot \Vert _{X}$ such that 
$\xi _{k}\in X$, $\forall k\in \mathbb{N}$. Clearly, $\Xi _{n}(X):=\mathrm{%
lin}\{\xi _{1},\cdots ,\xi _{n}\}\subset X$, $\forall n\in \mathbb{N}$ is a
sequence of closed subspaces of $X$ with the norm induced by $X$. Consider
the coordinate isomorphism $J$ defined as 
\[
\begin{array}{ccc}
J:\,\,\mathbb{R}^{n} & \longrightarrow & \Xi _{n}(X) \\ 
\alpha =(\alpha _{1},\cdots ,\alpha _{n}) & \longmapsto & 
\sum_{k=1}^{n}\alpha _{k}\cdot \xi _{k}.%
\end{array}%
\]%
Hence, the definition 
\[
\Vert \alpha \Vert _{J^{-1}\Xi _{n}(X)}=\Vert J\alpha \Vert _{X} 
\]%
induces the norm on $\mathbb{R}^{n}$. To be able to apply methods of
geometry of Banach spaces to various open problems in different spaces of
functions on $\mathbb{S}^{d}(\mathbb{C})$ we will need to calculate an
expectation of the function $\rho _{n}(\alpha ):=\Vert \alpha \Vert
_{J^{-1}\Xi _{n}(X)}$ on the unit sphere $\mathbb{S}^{n-1}\subset \mathbb{R}%
^{n}$ with respect to the invariant probabilistic measure $d\mu _{n}$, i.e.,
to find the Levy mean 
\[
M(\Vert \cdot \Vert _{J^{-1}\Xi _{n}(X)})\,=\,\int_{\mathbb{S}^{n-1}}\,\Vert
\alpha \Vert _{J^{-1}\Xi _{n}(X)}\cdot \,d\mu _{n}(\alpha ). 
\]%
As a motivating example consider the case $X=L_{p}:=L_{p}\left( \mathbb{S}%
^{d}(\mathbb{C})\right) $, 
\[
\Vert \phi \Vert _{p}:=\left\{ 
\begin{array}{cc}
\left( \int_{\mathbb{S}^{d}(\mathbb{C})}\left\vert \phi \right\vert
^{p}\cdot d\nu \right) ^{1/p}, & 1\leq p<\infty , \\ 
\mathrm{ess}\,\mathrm{sup}\,\left\vert \phi \right\vert , & p=\infty .%
\end{array}%
\right. 
\]%
In this case we shall write $\Vert \alpha \Vert _{(p)}=\Vert J\alpha \Vert
_{p}$. The sequence $\{\xi _{k}\}_{k\in \mathbb{N}}$ of orthonormal
harmonics on $\mathbb{S}^{d}(\mathbb{C})$ is not uniformly bounded on $%
\mathbb{S}^{d}(\mathbb{C})$. Hence, the method of estimating of Levy means
developed in \cite{ku1} - \cite{klw} can not give sharp order result.
Various modifications of this method presented in \cite{ku3} - \cite{ku5}
give an extra $(\log n)^{1/2}$ factor even if $p<\infty $. Our general
result concentrated in Lemma 3 which gives sharp order estimates for the
Levy means which correspond to the norm induced on $\mathbb{R}^{n}$ by the
subspace $\oplus _{s=1}^{m}\mathrm{H}_{k_{s}}\cap L_{p}$, $\mathrm{dim}%
\,\oplus _{s=1}^{m}\mathrm{H}_{k_{s}}=n$ with an arbitrary index set $%
(k_{1},\cdots ,k_{m})$, where $\mathrm{H}_{k_{s}}$ are the eigenspaces of
the Laplace-Beltrami operator for $\mathbb{S}^{d}(\mathbb{C})$ defined by (%
\ref{lb}). To show the boundness of the respective Levy means as $%
n\rightarrow \infty $ we employ the equality 
\[
\int_{\mathbb{R}^{n}}\,h(\alpha )\,d\gamma (\alpha )=\lim_{m\rightarrow
\infty }\,\,\int_{0}^{1}\,h\left( \frac{\delta _{1}^{m}(\theta )}{(2\pi
)^{1/2}},\cdots ,\frac{\delta _{n}^{m}(\theta )}{(2\pi )^{1/2}}\right)
d\theta , 
\]%
where $h:\mathbb{R}^{n}\rightarrow \mathbb{R}$ is a continuous function such
that 
\[
h(\alpha _{1},\cdots ,\alpha _{n})\,\mathrm{exp}\left(
-\sum_{k=1}^{n}|\alpha _{k}|\right) \,\rightarrow \,0 
\]%
uniformly when $\sum_{k=1}^{n}|\alpha _{k}|\,\rightarrow \infty $,

\[
d\gamma (\alpha )=\mathrm{exp}\left( -\pi \sum_{k=1}^{n}\alpha
_{k}^{2}\right) d\alpha 
\]%
is the Gaussian measure on $\mathbb{R}^{n}$,

\[
\delta _{k}^{m}(\theta ):=m^{-1/2}\cdot \left( r_{(k-1)m}(\theta )+\cdots
+r_{km}(\theta )\right) ,1\leq k\leq n 
\]%
and $r_{s}(\theta )=\mathrm{sign}\sin (2^{s}\pi \theta )$, $s\in \mathbb{N}%
\cup \{0\}$, $\theta \in \lbrack 0,1]$ is the sequence of Rademacher
functions \cite{kwap}, \cite{kuto1}. To extend our estimates to the case $%
p=\infty $ we apply Lemma 2 which gives a useful inequality between $1\leq
p,q\leq \infty $ norms of polynomials on $\mathbb{S}^{d}(\mathbb{C})$ with
an arbitrary spectrum. It seems that the factor $(\log n)^{1/2}$ obtained in
Lemma 3 is essential because of the lower bound for the Levy means found in 
\cite{kashin} in the case of trigonometric system. This fact explains a
logarithmic slot in our estimates presented in Theorem 2. We derive lower
bounds for $m-$therm approximation of Sobolev's classes (\ref{sobolev})
using Lemmas 1 and 2, Urysohn's inequality, Bourgain-Milman inequality and
estimates of Levy means given by Lemma 3 and (\ref{infty}). Upper bounds for 
$m-$therm approximation contained in Theorem 1 where we establish Jackson's
type inequality. As it follows from Remark 1, $m $-term approximations can
not give better rate of convergence than approximation by the subspace of
polynomials on $\mathbb{S}^{d}(\mathbb{C})$ of the same dimension.

In this article there are several universal constants which enter into the
estimates. These positive constants are mostly denoted by $C,C_{1},...$. We
will only distinguish between the different constants where confusion is
likely to arise, but we have not attempted to obtain good estimates for
them. For ease of notation we will write $a_{n}\ll b_{n}$ for two sequences,
if $a_{n}\leq C\cdot b_{n}$, $\forall n\in \mathbb{N}$ and $a_{n}\asymp
b_{n} $, if $C_{1}\cdot b_{n}\leq a_{n}\leq C_{2}\cdot b_{n}$, $\forall n\in 
\mathbb{N}$ and some constants $C$, $C_{1}$ and $C_{2}$. Also, we shall put $%
(a)_{+}:=\max \{a,0\}$.


\section{Harmonic Analysis}

\label{harmonic analysis}

Let $\mathbb{C}^{n}$ be $n$-dimensional complex space. We will denote
vectors in $\mathbb{C}^{n}$ by $\mathbf{z}=(z_{1},...,z_{n})$. Let the inner
product of two vectors $\mathbf{w,z}\in \mathbb{C}^{n}$ be 
\[
\langle \mathbf{w},\mathbf{z}\rangle =\sum_{j+1}^{n}w_{j}\cdot \overline{z}%
_{j}, 
\]%
and the length of a vector be $|\mathbf{z}|=\langle \mathbf{w},\mathbf{z}%
\rangle ^{1/2}$. Let 
\[
\mathbb{S}^{d}(\mathbb{C}):=\{\mathbf{z}\in \mathbb{C}^{n}:\,|\mathbf{z}%
|=1\} 
\]%
be the unit sphere in $\mathbb{C}^{n}$. Here $d$ means the topological
dimension of the complex sphere over reals. It means that $d=2n-1$.

Observe that $\mathbb{S}^{d}(\mathbb{C})$ is a compact, connected, $d$%
-dimensional, $C^{\infty }$ Riemannian manifold with $C^{\infty }$ metric.
Let $g$ its metric tensor, $\nu $ its normalized volume element and $\Delta $
its Laplace-Beltrami operator. In local coordinates $x_{l}$, $1\leq l\leq d$%
, 
\begin{equation}
\Delta =-(\overline{g})^{-1/2}\cdot \sum_{k}\frac{\partial }{\partial x_{k}}%
\left( \sum_{j}g^{jk}\cdot (\overline{g})^{1/2}\cdot \frac{\partial }{%
\partial x_{j}}\right) ,  \label{lb}
\end{equation}%
where $g_{jk}:=g(\partial /x_{j},\partial /x_{k})$, $\overline{g}:=|\mathrm{%
det}(g_{jk})|$, and $(g^{jk}):=(g_{jk})^{-1}$. It is well-known that $\Delta 
$ is an elliptic, self adjoint, invariant under isometry, second order
operator. The eigenvalues $\theta _{k}=k\cdot \left( k+d-1\right) $, of $%
\Delta $ are discrete, nonnegative and form an increasing sequence $0\leq
\theta _{0}\leq \theta _{1}\leq \cdots \leq \theta _{n}\leq \cdots $ with $%
+\infty $ the only accumulation point. The corresponding eigenspaces $%
\mathrm{H}_{k}$, $k\geq 0$ are finite-dimensional, orthogonal with respect
to the scalar product 
\[
\langle f,g\rangle :=\int_{\mathbb{S}^{d}(\mathbb{C})}f\,\,\cdot \overline{g}%
\cdot d\nu 
\]%
and 
\[
L_{2}:=L_{2}(\mathbb{S}^{d}(\mathbb{C}),\nu )=\mathrm{cl}_{L_{2}}\left(
\bigoplus_{k=0}^{\infty }\mathrm{H}_{k}\right) . 
\]%
It is known \cite{lk} that $d_{n}:=\mathrm{dim}\,(\mathrm{H}_{n})\asymp
n^{2d-1}.$

The complex sphere $\mathbb{S}^{d}(\mathbb{C})$, $d=3,5,...$ is invariant
under the action of the unitary group $\mathcal{U}_{(d+1)/2},$ the group of $%
(d+1)/2\times (d+1)/2$ complex matrices $U$ which satisfy $UU^{\ast
}=I_{(d+1)/2},$ where $U_{ij}^{\ast }=\overline{U_{ji}}$, $1\leq i,j\leq
(d+1)/2$ and $\mathbb{S}^{d}(\mathbb{C})=\mathcal{U}_{(d+1)/2}/\mathcal{U}%
_{(d-1)/2}.$ If $\kappa $ is a $\mathcal{U}_{(d+1)/2}$ invariant kernel then
there is a univariate function $\Psi $ such that $\kappa (\mathbf{x},\mathbf{%
y})=\Psi (\left\langle \mathbf{x},\mathbf{y}\right\rangle ).$ We define the
convolution of $f\in L_{1}\left( \mathbb{S}^{d}(\mathbb{C})\right) $ with a $%
\mathcal{U}_{(d+1)/2}$-invariant kernel $\kappa $ as 
\[
(f\ast \kappa )(\mathbf{x})=\int_{\mathbb{S}^{d}(\mathbb{C})}f(\mathbf{y}%
)\cdot \Psi (\left\langle \mathbf{x},\mathbf{y}\right\rangle )\cdot d\nu
(y). 
\]%
Let $M_{k}$ be an invariant kernel of orthogonal projector $L_{2}\rightarrow 
\mathrm{H}_{k}$. Then $M_{k}\ast \phi \in \mathrm{H}_{k}$ for any $\phi \in
L_{2}$. Let us fix an orthonormal basis $\{Y_{m}^{k}\}_{m=1}^{d_{k}}$ of $%
\mathrm{H}_{k}$. For an arbitrary $\phi \in L_{p}$, $1\leq p\leq \infty $
with the formal Fourier series 
\[
\phi \sim \sum_{k\in \mathbb{N}\cup \{0\}}M_{k}\ast \phi =\sum_{k\in \mathbb{%
N}\cup \{0\}}\sum_{m=1}^{d_{k}}c_{k,m}(\phi )\cdot
Y_{m}^{k},\,\,\,c_{k,m}(\phi )=\int_{\mathbb{S}^{d}(\mathbb{C})}\phi \cdot 
\overline{Y_{m}^{k}}d\nu , 
\]%
the $\gamma $-th fractional integral $I_{\gamma }\phi :=\phi _{\gamma }$, $%
\gamma >0$, is defined as 
\begin{equation}
\phi _{\gamma }\sim C+\sum_{k\in \mathbb{N}}\theta _{k}^{-\gamma
/2}\sum_{m=1}^{d_{k}}c_{k,m}(\phi )\cdot Y_{m}^{k},\,\,\,C\in \mathbb{R}%
.\,\,\,\   \label{sobolev}
\end{equation}%
The function $D_{\gamma }\phi :=\phi ^{(\gamma )}\in L_{p}$, $1\leq p\leq
\infty $ is called the $\gamma $-th fractional derivative of $\phi $ if 
\[
\phi ^{(\gamma )}\sim \sum_{k\in \mathbb{N}}\theta _{k}^{\gamma
/2}\sum_{m=1}^{d_{k}}c_{k,m}(\phi )\cdot Y_{m}^{k}. 
\]%
The Sobolev classes $W_{p}^{\gamma }$ are defined as sets of functions with
formal Fourier expansions (\ref{sobolev}) where $\Vert \phi \Vert _{p}\leq 1$
and $\int_{\mathbb{M}^{d}}\phi d\nu =0$.

We recall that a Riemannian manifold $\mathbb{M}^{d}$ is called homogeneous
if its group of isometries ${\mathcal{G}}$ acts transitively on it, i.e. for
every $x,y\in \mathbb{M}^{d}$, there is a $g\in {\mathcal{G}}$ such that $%
gx=y$. For a compact homogeneous Riemannian manifold $\mathbb{M}^{d}$ which
is, in particular, $\mathbb{S}^{d}(\mathbb{C})$ the following addition
formula is known \cite{gine} 
\begin{equation}
\sum_{k=1}^{d_{k}}|Y_{m}^{k}(x)|^{2}=d_{k},\,\,\,\forall x\in \mathbb{M}^{d},
\label{addi}
\end{equation}%
where $\{Y_{m}^{k}\}_{m=1}^{d_{k}}$ is an arbitrary orthonormal basis of $%
\mathrm{H}_{k}$, $k\geq 0$.

\section{$m$-Term Approximation}

Our upper bounds come from Jackson's type inequality.

\textbf{Theorem 1} \emph{\ Let $f\in L_{p}$ and 
\[
E(f,\mathcal{T}_{N},L_{p}):=\inf_{t_{N}\in \mathcal{T}_{N}}\,\,\Vert
f-t_{N}\Vert _{p}. 
\]%
be the best approximation of $f$ by $\mathcal{T}_{N}$. If $f^{(\gamma )}\in
L_{p}$ and $\gamma >(d-1)/2$ then 
\[
E(f,\mathcal{T}_{N},L_{p})\leq C\cdot N^{-\gamma }\cdot E(f^{(\gamma )},%
\mathcal{T}_{N},L_{p}),\,\,\,1\leq p\leq \infty . 
\]%
} \textbf{Proof} To produce our estimates we will need some information
concerning Ces\`{a}ro means. The Ces\`{a}ro kernel is defined by 
\[
S_{n}^{\delta }:=\frac{1}{C_{n}^{\delta }}\sum_{m=0}^{n}C_{n-m}^{\delta
}\cdot M_{m}, 
\]%
where $C_{n}^{\delta }$ are Ces\`{a}ro numbers of order $n$ and index $%
\delta $, i.e. 
\begin{equation}
C_{n}^{\delta }=\frac{\Gamma (n+\delta +1)}{\Gamma (\delta +1)\cdot \Gamma
(n+1)}\asymp n^{\delta }.  \label{delta}
\end{equation}%
It is known \cite{lk} that for $0\leq \delta \leq (d+1)/2,$ 
\begin{equation}
\left\Vert S_{n}^{\delta }\right\Vert _{1}\leq C\left\{ 
\begin{array}{cc}
n^{(d-1)/2-\delta }, & \delta \leq (d-3)/2, \\ 
\left( \log n\right) ^{2}, & \delta =(d-1)/2, \\ 
1, & \delta =(d+1)/2.%
\end{array}%
\right.  \label{cesaro}
\end{equation}%
Fix a polynomial $\phi _{M}\in \mathcal{T}_{M}$ with $\Vert \phi _{M}\Vert
_{p}\leq 1$ and let 
\[
K_{N}:=\sum_{k=1}^{N}\lambda _{k}\cdot M_{k}. 
\]%
Let $\left\{ \lambda _{k}\right\} _{k\in \mathbb{N}}$ be a fixed sequence of
real numbers. Applying Abel's transform $s+1$ times where $s:=(d+1)/2$ we
see that, for $N>s+1$, 
\[
K_{N}\ast \phi _{M}=(K_{N-s-1}^{1}+K_{N}^{2})\ast \phi _{M}, 
\]%
where 
\[
K_{N}^{1}:=\sum_{k=1}^{N}\Delta ^{s+1}\lambda _{k}\cdot C_{k}^{s}\cdot
S_{k}^{s}, 
\]%
\[
K_{N}^{2}:=\sum_{k=0}^{(d+1)/2}\Delta ^{k}\lambda _{N-k}\cdot
C_{N-k}^{k}\cdot S_{N-k}^{k}, 
\]%
$\Delta ^{0}\lambda _{k}:=\lambda _{k}$, $\Delta ^{1}\lambda _{k}=\lambda
_{k}-\lambda _{k+1}$ and $\Delta ^{s+1}\lambda _{k}=\Delta ^{s}\lambda
_{k}-\Delta ^{s}\lambda _{k+1}$, $k,s\in \mathbb{N}$. Using (\ref{delta})
and (\ref{cesaro}) we get 
\[
\left\Vert K_{N}^{1}\right\Vert _{1}\leq \sum_{k=1}^{N}\left\vert \Delta
^{s+1}\lambda _{k}\right\vert \cdot C_{k}^{s}\cdot \left\Vert
S_{K}^{s}\right\Vert _{1}\leq C\cdot \sum_{k=1}^{N}\left\vert \Delta
^{s+1}\lambda _{k}\right\vert \cdot k^{s}\cdot \left\Vert
S_{k}^{s}\right\Vert _{1} 
\]%
\[
\leq C\cdot \sum_{k=1}^{N}\left\vert \Delta ^{s+1}\lambda _{k}\right\vert
\cdot k^{s}=C\cdot \sum_{k=1}^{N}\left\vert \Delta ^{(d+3)/2}\lambda
_{k}\right\vert \cdot k^{(d+1)/2} 
\]%
\begin{equation}
\leq C\cdot \sum_{k=1}^{N}k^{-\gamma -(d+3)/2}\cdot
k^{(d+1)/2}=\sum_{k=1}^{N}k^{-\gamma -1}\leq C,\,\,\,\gamma >0.  \label{k1}
\end{equation}%
Since in our case $\lambda _{k}=\theta _{k}^{-\gamma /2}$ then $|\Delta
^{(d+3)/2}\lambda _{k}|\asymp k^{-\gamma -(d+3)/2}$ as $k\rightarrow \infty $
and, by (\ref{cesaro}), $\left\Vert S_{k}^{s}\right\Vert _{1}\leq C$ as $%
k\rightarrow \infty $ . Similarly, using (\ref{cesaro}) we get 
\[
\Vert K_{N}^{2}\Vert _{1}\leq \sum_{k=0}^{(d+1)/2}|\Delta ^{k}\lambda
_{N-k}|\cdot C_{N-k}^{k}\cdot \Vert S_{N-k}^{k}\Vert _{1} 
\]%
\[
\leq C\cdot \sum_{k=0}^{(d+1)/2}|\Delta ^{k}\lambda _{N-k}|\cdot
(N-k)^{k}\cdot (N-k)^{(d+1)/2-1-k} 
\]%
\begin{equation}
\leq C\cdot N^{-\gamma +(d-1)/2}.  \label{k2}
\end{equation}%
From (\ref{k2}) follows that if $\gamma >(d-1)/2$ then 
\begin{equation}
\lim_{n\rightarrow \infty }\Vert K_{N}^{2}\Vert _{1}=0.  \label{k3}
\end{equation}%
Comparing (\ref{k1}) and (\ref{k3}) we get that for any fixed polynomial $%
\phi _{M}\in \mathcal{T}_{M}$, $M\in \mathbb{N}$ the sequence of functions $%
K_{N}^{1}\ast \phi _{M}$ converges in $L_{1}$ to the function 
\[
K\ast \phi _{M}=\left( \sum_{k=1}^{\infty }\theta _{k}^{-\gamma /2}\cdot
M_{k}\right) \ast \phi _{M}. 
\]%
Remark that 
\begin{equation}
\left\Vert K-K_{N}^{1}\right\Vert _{1}\leq \sum_{k=N+1}^{\infty }\left\vert
\Delta ^{s+1}\lambda _{k}\right\vert \cdot C_{k}^{s}\cdot \left\Vert
S_{k}^{s}\right\Vert _{1}\leq CN^{-\gamma }.  \label{gam}
\end{equation}%
Fix an arbitrary polynomial $\psi _{N}\in \mathcal{T}_{N}.$ For any $f$, $%
f=K\ast f^{(\gamma )}$ such that $f^{(\gamma )}\in L_{p}$ we have 
\[
E(f,\mathcal{T}_{N},L_{p})\leq \left\Vert K\ast f^{(\gamma )}-K\ast \psi
_{N}+K_{N}^{1}\ast \psi _{N}-K_{N}^{1}\ast f^{(\gamma )}\right\Vert _{p} 
\]%
\[
=\left\Vert K\ast \left( f^{(\gamma )}-\psi _{N}\right) -K_{N}^{1}\ast
\left( f^{(\gamma )}-\psi _{N}\right) \right\Vert _{p} 
\]%
\[
\leq \left\Vert (K-K_{N}^{1})\ast \left( f^{(\gamma )}-\psi _{N}\right)
\right\Vert _{p}\leq \left\Vert (K-K_{N}^{1})\right\Vert _{1}\cdot
\left\Vert \left( f^{(\gamma )}-\psi _{N}\right) \right\Vert _{p} 
\]%
\[
\leq C\cdot N^{-\gamma }\cdot E(f^{(\gamma )},\mathcal{T}_{N},L_{p}), 
\]%
where in the last line we used (\ref{gam}) and the fact that $\psi _{N}$ is
an arbitrary polynomial. $\blacksquare $

\textbf{Remark 1} \emph{From Theorem 1, (\ref{ber}) and \cite{dz}, p. 658 we
get Kolmogorov's type inequality, 
\[
\Vert f^{(\alpha )}\Vert _{p}\leq C\Vert f^{(\beta )}\Vert _{p}^{\alpha
/\beta }\cdot \Vert f\Vert _{p}^{1-\alpha /\beta }, 
\]%
where $1\leq p\leq \infty $ and $(d-1)/2\leq \alpha \leq \beta $. }

To prove our lower bounds we will need several Lemmas.

\textbf{Lemma 1} \emph{There is a sequence of function $Q_{2N}\in \mathcal{T}%
_{2N},N\in \mathbb{N}$ such that for any $t_{N}\in \mathcal{T}_{N}$ we have 
\[
Q_{2N}\ast t_{N}=t_{N} 
\]%
and 
\[
\Vert Q_{2N}\Vert _{1}\leq C,\,\,\,\forall N\in \mathbb{N}. 
\]%
}

The proof of this statement is based on the norm estimates for the Ces\`{a}%
ro means (\ref{cesaro}) and the line of arguments used in \cite{klt2009}
(see Example 2.6). Observe that $Q_{2N}$ has the form 
\[
Q_{2N}=\sum_{k=0}^{2N}\lambda _{k}^{\left( 2N\right) }\cdot M_{k}, 
\]%
where $M_{k}$ is the reproducing kernel for $\mathrm{H}_{k}$, 
\[
\lambda _{k}^{\left( 2N\right) }=\chi _{d}\left( \frac{k}{2N}\right) ,0\leq
k\leq 2N, 
\]%
\[
\chi _{s}\left( t\right) =2d\cdot \int_{t}^{t+1/\left( 2d\right) }\chi
_{s-1}\left( u\right) \cdot du,1\leq s\leq d, 
\]%
and 
\[
\chi _{0}\left( t\right) =\left\{ 
\begin{array}{cc}
1, & t\in \left[ 0,1\right] , \\ 
0, & t\notin \left[ 0,1\right] .%
\end{array}%
\right. 
\]%
The function $\chi _{d}\left( t\right) $ is $d-1$ times continuously
differentiable and nonnegative on $\left[ 0,\infty \right) .$ Furthermore, $%
\chi _{d}^{\left( d-1\right) }\left( t\right) $ is Lipschitz continuous, $%
\chi _{d}\left( t\right) =1$ if $t\in \left[ 0\leq t\leq 1/2\right] ,$ and 
\[
\chi _{d}\left( t\right) =\frac{\left( 2d\right) ^{d}}{d!}\left( 1-t\right)
^{d},1-\frac{1}{2d}\leq t\leq 1. 
\]%
Also $\chi _{d}\left( t\right) $ is a polynomial of degree $d$ in each
interval $\left[ t_{s},t_{s-1}\right] ,1\leq s\leq d,$ where $%
t_{s}=1-s/\left( 2d\right) .$

Let $\Omega _{m}:=\{k_{1}<\cdots <k_{m}\}\subset \mathbb{N}$ and $\Xi
_{n}(\Omega _{m}):=\mathrm{lin}\{\mathrm{H}_{k_{l}}\}_{l=1}^{m}.$

\textbf{Lemma 2} \emph{For any $\Omega _{m}$ and any $\xi \in \Xi
_{n}(\Omega _{m})$, $m\in \mathbb{N}$ we have 
\[
\Vert \xi \Vert _{q}\leq n^{(1/p-1/q)_{+}}\cdot \Vert \xi \Vert _{p},\,\,\, 
\]%
where $1\leq p,q\leq \infty $ and $n:=\mathrm{dim}\,\Xi _{n}(\Omega _{m})$. }

\textbf{Proof} Let 
\[
K_{n}(x,y):=\sum_{i=1}^{n}\eta _{i}(x)\overline{\cdot \eta _{i}(y)}. 
\]%
be the reproducing kernel for $\Xi _{n}(\Omega _{m})$. Clearly, 
\[
K_{n}(x,y)=\int_{\mathbb{S}^{d}(\mathbb{C})}K_{n}(x,z)\cdot K_{n}(z,y)\cdot
d\nu (z), 
\]%
and $K_{n}(x,y)=\overline{K_{n}(y,x)}$. Hence, using the Cauchy-Schwartz
inequality, 
\[
\Vert K_{n}(\cdot ,\cdot )\Vert _{\infty }\leq \Vert K_{n}(y,\cdot )\Vert
_{2}\cdot \Vert K_{n}(x,\cdot )\Vert _{2} 
\]%
for any $x,y\in \mathbb{S}^{d}(\mathbb{C})$. Due to the addition formula (%
\ref{addi}), we have $\Vert K_{n}(x,\cdot )\Vert _{2}=n^{1/2}$. Therefore, 
\begin{equation}
\Vert K_{n}(\cdot ,\cdot )\Vert _{\infty }\leq n.  \label{2010}
\end{equation}%
Let $\xi \in \Xi (\Omega _{m})$. Then applying H\"{o}lder inequality and (%
\ref{2010}) we get 
\[
\Vert \xi \Vert _{\infty }\leq \Vert K_{n}(\cdot ,\cdot )\Vert _{\infty
}\cdot \Vert \xi \Vert _{1}\leq n\cdot \Vert \xi \Vert _{1}, 
\]%
and hence 
\[
\Vert I\Vert _{L_{1}(\mathbb{M}^{d})\cap \Xi _{n}(\Omega _{m})\rightarrow
L_{\infty }(\mathbb{M}^{d})\cap \Xi _{n}(\Omega _{m})}\leq n, 
\]%
where $I:L_{p}\rightarrow L_{q}$ is the embedding operator. Trivially, 
\[
\Vert I\Vert _{L_{p}(\mathbb{M}^{d})\cap \Xi (\Omega _{m})\rightarrow L_{p}(%
\mathbb{M}^{d})\cap \Xi (\Omega _{m})}=1, 
\]%
where $1\leq p\leq \infty $. Hence, using the Riesz-Thorin interpolation
Theorem and embedding arguments we obtain 
\[
\Vert \xi \Vert _{p}\leq n^{(1/p-1/q)_{+}}\cdot \Vert \xi \Vert
_{q},\,\,\,\forall \xi \in \Xi _{n}(\Omega _{m}),\,\,1\leq p,q\leq \infty
.\,\, 
\]

$\blacksquare $

Let $\mathbb{R}^{n}$ be the Euclidean space with the usual scalar product $%
\left\langle \alpha ,\beta \right\rangle :=\sum_{k=1}^{n}\alpha _{k}\beta
_{k},$ and the norm $\left\vert \alpha \right\vert :=\left\langle \alpha
,\alpha \right\rangle ^{1/2}.$ Let $B_{(2)}^{n}:=\left\{ \alpha \left\vert
\alpha \in \mathbb{R}^{n},\left\vert \alpha \right\vert \leq 1\right.
\right\} $ be the canonical Euclidean ball and $\mathbb{S}^{n-1}:=\left\{
\alpha \left\vert \alpha \in \mathbb{R}^{n},\left\vert \alpha \right\vert
=1\right. \right\} .$ Let $K\subset $ $\mathbb{R}^{n},$ be a convex,
centrally symmetric body, then the polar set $K^{o}$ is defined as $%
K^{o}:=\sup \left\{ \left\vert \left\langle \alpha ,\beta \right\rangle
\right\vert \leq 1,\beta \in K\right\} .$ Let us fix a norm $\Vert \cdot
\Vert $ on $\mathbb{R}^{n}$ and let $E=(\mathbb{R}^{n},\Vert \cdot \Vert )$
with the unit ball $B_{E}$. The dual space $E^{o}=(\mathbb{R}^{n},\Vert
\cdot \Vert ^{o})$ is endowed with the norm $\Vert \xi \Vert
^{o}=\sup_{\sigma \in B_{E}}|\langle \xi ,\sigma \rangle |$ and has the unit
ball $B_{E^{o}}:=\left( B_{E}\right) ^{o}$. In these notations the Levy mean 
$M_{B_{E}}$ is 
\[
M_{B_{E}}=\int_{\mathbb{S}^{n-1}}\Vert \xi \Vert d\mu _{n}, 
\]%
where $d\mu _{n}$ denotes the normalized invariant measure on $\mathbb{S}%
^{n-1}$, the unit sphere in $\mathbb{R}^{n}$. We are interested in the case
where $\Vert \cdot \Vert =\left\Vert \cdot \right\Vert _{(p)}$. In this case
we shall write $J^{-1}B_{L_{p}\cap \Xi _{n}(\Omega _{m})}=B_{(p)}^{n}$ In
the case $\Omega _{m}=\{1,\cdots ,m\}$ the estimates of the associated Levy
means were obtained in \cite{kuto1}. This result can be easily generalized
to an arbitrary index set $\Omega _{m}=\{k_{1}<\cdots <k_{m}\}.$

\textbf{Lemma 3} \emph{%
\[
M_{B_{(p)}^{n}}\leq C\cdot p^{1/2},\,\,p<\infty . 
\]%
}

Applying Lemmas 2 and 3 with $p=\log n$ we get 
\[
M_{B_{(\infty )}^{n}}=\int_{\mathbb{S}^{n-1}}\Vert \xi \Vert _{(\infty
)}\cdot d\mu _{n}\leq n^{1/p}\cdot \int_{\mathbb{S}^{n-1}}\Vert \xi \Vert
_{(p)}\cdot d\mu _{n} 
\]%
\begin{equation}
\leq C\cdot p^{1/2}\cdot n^{1/p}=C\cdot \left( \log n\right) ^{1/2}\cdot
n^{1/(\log n)}\leq C\cdot \left( \log n\right) ^{1/2}.  \label{infty}
\end{equation}%
Our lower bounds for $m$-term approximation are concentrated in

\textbf{Theorem 2 } \emph{%
\[
\nu _{m}(W_{p}^{\gamma },\Xi ,L_{q})\geq C\cdot m^{-\gamma /d}\cdot
\vartheta _{m},
\]%
where 
\[
\vartheta _{m}\geq C\left\{ 
\begin{array}{cc}
p^{1/2}\cdot 2^{C(q{\prime })^{1/2}}, & 1<q,p<\infty , \\ 
(\log m)^{1/2}\cdot 2^{C(q^{\prime })^{1/2}}, & p=\infty ,q>1,%
\end{array}%
\right\} ^{-C},1/q+1/q^{\prime }=1.
\]%
}

\textbf{Proof} It is sufficient to consider the case $p\geq 2$ and $1\leq
q\leq 2$ since all other cases follow by embedding arguments. By Bernstein's
inequality (\ref{ber}), 
\[
(\dim \mathcal{T}_{N})^{-\gamma /d}U_{p}\cap \mathcal{T}_{N}\subset
W_{p}^{\gamma }. 
\]%
Hence, from the definition of $\nu _{m}$ it follows 
\[
\nu _{m}(W_{p}^{\gamma },\Xi ,L_{q})\geq \nu _{m}((\dim \mathcal{T}%
_{N})^{-\gamma /d}\cdot U_{p}\cap \mathcal{T}_{N},\Xi ,L_{q}) 
\]%
\begin{equation}
=(\dim \mathcal{T}_{N})^{-\gamma /d}\cdot \nu _{m}(U_{p}\cap \mathcal{T}%
_{N},\Xi ,L_{q}).  \label{x}
\end{equation}%
Let $\phi \in U_{p}\cap \mathcal{T}_{N}$ and $\xi \in L_{p}$. Then applying
Lemma 1 we get 
\[
\Vert Q_{2N}\ast (\phi -\xi )\Vert _{p}=\Vert Q_{2N}\ast \phi -Q_{2N}\ast
\xi \Vert _{p}\leq \Vert Q_{2N}\Vert _{1}\cdot \Vert \phi -\eta \Vert _{p}, 
\]%
where $\eta :=Q_{2N}\ast \xi \in \mathcal{T}_{2N}$ and $\phi =Q_{2N}\ast
\phi $ for any $\phi \in \mathcal{T}_{N}$. Consequently, 
\begin{equation}
\inf_{\eta \in \Xi (\Omega _{m})}\,\Vert \phi -\eta \Vert _{p}\geq \frac{1}{%
\Vert Q_{2N}\Vert _{1}}\,\inf_{\xi \in \Xi (\Omega _{m})\cap \mathcal{T}%
_{2N}}\Vert \phi -\xi \Vert _{p}\geq C\inf_{\eta \in \Xi (\Omega _{m})\cap 
\mathcal{T}_{2N}}\Vert \phi -\eta \Vert _{p},  \label{xx}
\end{equation}%
for any $\phi \in U_{p}\cap \mathcal{T}_{N}$. Comparing (\ref{x}) and (\ref%
{xx}) we find 
\begin{equation}  \label{bumbum}
\nu _{m}(W_{p}^{\gamma },\Xi ,L_{q}(\mathbb{M}^{d}))\geq (\dim \mathcal{T}%
_{N})^{-\gamma /d}\cdot \vartheta _{m},  \label{00}
\end{equation}%
where 
\[
\vartheta _{m}:=\nu _{m}(U_{p}\cap \mathcal{T}_{N},\Xi \cap \mathcal{T}%
_{2N},L_{q}\cap \mathcal{T}_{2N}). 
\]%
Let $e_{1},\cdots ,e_{k}$ be the canonic basis in $\mathbb{R}^{k}$, $k:=\dim 
\mathcal{T}_{2N}$. Let $\mathbf{l}=\{k_{1},\cdots ,k_{m}\}\in \mathbb{N}^{m}$%
, $k_{s}\leq k$, $1\leq s\leq m$ and $X_{\mathbf{l}}^{m}=\mathrm{lin}%
\{e_{k_{s}}\}_{s=1}^{m}$. Since $p\geq 2,$ then by H\"{o}lder's inequality $%
B_{p}\subset B_{2}.$ Consequently, $J^{-1}(B_{p}\cap \mathcal{T}_{N})\subset
J^{-1}(B_{2}\cap \mathcal{T}_{N})$ and, therefore, 
\[
J^{-1}(B_{p}\cap \mathcal{T}_{N})\subset \cup _{\mathbf{l}}(X_{\mathbf{l}%
}^{m}+\,J^{-1}(\left( \vartheta _{m}\cdot B_{q}\right) \,\cap \mathcal{T}%
_{2N}))\,\cap J^{-1}(B_{p}\cap \mathcal{T}_{N}) 
\]%
\begin{equation}
\subset \cup _{\mathbf{l}}(X_{\mathbf{l}}^{m}+\,J^{-1}(\left( \vartheta
_{m}\cdot B_{q}\right) \cap \mathcal{T}_{2N}))\cap J^{-1}(B_{2}\cap \mathcal{%
T}_{N}).  \label{b}
\end{equation}%
Let $P\left( X_{\mathbf{l}}^{m}\right) $ be the orthoprojector onto $X_{%
\mathbf{l}}^{m}$ and $P^{\perp }\left( X_{\mathbf{l}}^{m}\right) $ be the
orthoprojector onto $\left( X_{\mathbf{l}}^{m}\right) ^{\perp }$ in $J^{-1}%
\mathcal{T}_{N}$. Observe that for any $\mathbf{l}$ 
\[
(X_{\mathbf{l}}^{m}+\,J^{-1}(\left( \vartheta _{m}\cdot B_{q}\right) \cap 
\mathcal{T}_{2N}))\cap J^{-1}(B_{2}\cap \mathcal{T}_{N}) 
\]%
\[
\subset P\left( X_{\mathbf{l}}^{m}\right) \left( X_{\mathbf{l}}^{m}\cap
J^{-1}(B_{2}\cap \mathcal{T}_{N})\right) +\,P^{\perp }\left( X_{\mathbf{l}%
}^{m}\right) \circ J^{-1}(\left( \left( \vartheta _{m}\cdot B_{q}\right)
\cap \mathcal{T}_{N}\right) \cap \left( B_{2}\cap \mathcal{T}_{N}\right) ) 
\]%
\[
=X_{\mathbf{l}}^{m}\cap J^{-1}(B_{2}\cap \mathcal{T}_{N})+\,P^{\perp }\left(
X_{\mathbf{l}}^{m}\right) \circ J^{-1}(\left( \left( \vartheta _{m}\cdot
B_{q}\right) \cap \mathcal{T}_{N}\right) \cap \left( B_{2}\cap \mathcal{T}%
_{N}\right) ). 
\]%
Let 
\[
\dim \mathcal{T}_{N}:=n,l:=\dim J(X_{\mathbf{l}}^{m}\cap J^{-1}\mathcal{T}%
_{N})_{\mathcal{T}_{N}}^{\perp },\,\,\,s:=n-l. 
\]%
Taking volumes we get 
\[
\mathrm{Vol}_{n}((X_{\mathbf{l}}^{m}+\,J^{-1}(\left( \vartheta _{m}\cdot
B_{q}\right) \cap \mathcal{T}_{2N}))\cap J^{-1}(B_{2}\cap \mathcal{T}_{N})) 
\]%
\[
\leq \mathrm{Vol}_{n}\left( X_{\mathbf{l}}^{m}\cap J^{-1}(B_{2}\cap \mathcal{%
T}_{N})+\,P^{\perp }\left( X_{\mathbf{l}}^{m}\right) \circ J^{-1}(\left(
\left( \vartheta _{m}\cdot B_{q}\right) \cap \mathcal{T}_{N}\right) \cap
\left( B_{2}\cap \mathcal{T}_{N}\right) )\right) 
\]%
\[
=\mathrm{Vol}_{s}\,\left( B_{(2)}^{s}\right) \cdot \mathrm{Vol}_{l}\left(
P^{\perp }\left( X_{\mathbf{l}}^{m}\right) \circ J^{-1}(\left( \left(
\vartheta _{m}\cdot B_{q}\right) \cap \mathcal{T}_{N}\right) \cap \left(
B_{2}\cap \mathcal{T}_{N}\right) )\right) . 
\]%
To get an upper bound for 
\[
\mathrm{Vol}_{l}\left( P^{\perp }\left( X_{\mathbf{l}}^{m}\right) \circ
J^{-1}(\left( \left( \vartheta _{m}\cdot B_{q}\right) \cap \mathcal{T}%
_{N}\right) \cap \left( B_{2}\cap \mathcal{T}_{N}\right) )\right) 
\]%
we proceed as following. Let $x_{1},\cdot \cdot \cdot ,x_{N}$ be a $1$-net
for $J^{-1}\left( B_{q}\cap \mathcal{T}_{N}\right) $ in the norm induced by $%
J^{-1}\left( B_{2}\cap \mathcal{T}_{N}\right) .$ Hence, we have 
\[
J^{-1}\left( B_{q}\cap \mathcal{T}_{N}\right) \subset
\dbigcup\limits_{k=1}^{N}\left( x_{k}+J^{-1}\left( B_{2}\cap \mathcal{T}%
_{N}\right) \right) 
\]%
and, therefore, 
\[
\mathrm{Vol}_{l}\left( P^{\perp }\left( X_{\mathbf{l}}^{m}\right) \circ
J^{-1}(\left( \vartheta _{m}\cdot B_{q}\cap \mathcal{T}_{N}\right) \cap
\left( B_{2}\cap \mathcal{T}_{N}\right) )\right) 
\]%
\[
\leq \mathrm{Vol}_{l}\left( P^{\perp }\left( X_{\mathbf{l}}^{m}\right) \circ
J^{-1}(\vartheta _{m}\cdot B_{q}\cap \mathcal{T}_{N}\right) 
\]%
\[
=\vartheta _{m}^{l}\cdot \mathrm{Vol}_{l}\left( P^{\perp }\left( X_{\mathbf{l%
}}^{m}\right) \circ J^{-1}(B_{q}\cap \mathcal{T}_{N}\right) 
\]
\[
\leq \vartheta _{m}^{l}\cdot \mathrm{Vol}_{l}P^{\perp }\left( X_{\mathbf{l}%
}^{m}\right) \left( \dbigcup\limits_{k=1}^{N}\left( x_{k}+J^{-1}\left(
B_{2}\cap \mathcal{T}_{N}\right) \right) \right) 
\]%
\[
=\vartheta _{m}^{l}\cdot \mathrm{Vol}_{l}\left(
\dbigcup\limits_{k=1}^{N}P^{\perp }\left( X_{\mathbf{l}}^{m}\right) \left(
x_{k}+J^{-1}\left( B_{2}\cap \mathcal{T}_{N}\right) \right) \right) 
\]%
\[
\leq \vartheta _{m}^{l}\cdot \mathrm{Vol}_{l}\left(
\dbigcup\limits_{k=1}^{N}\left( P^{\perp }\left( X_{\mathbf{l}}^{m}\right)
x_{k}+P^{\perp }\left( X_{\mathbf{l}}^{m}\right) \circ J^{-1}\left(
B_{2}\cap \mathcal{T}_{N}\right) \right) \right) 
\]%
\[
\leq \vartheta _{m}^{l}\cdot \sum_{k=1}^{N}\mathrm{Vol}_{l}\left( P^{\perp
}\left( X_{\mathbf{l}}^{m}\right) \circ J^{-1}\left( B_{2}\cap \mathcal{T}%
_{N}\right) \right) 
\]

\[
=\vartheta _{m}^{l}\cdot N\cdot \mathrm{Vol}_{l}\left( B_{(2)}^{l}\right) ,
\]%
since $P^{\perp }\left( X_{\mathbf{l}}^{m}\right) \circ J^{-1}\left(
B_{2}\cap \mathcal{T}_{N}\right) =B_{(2)}^{l}.$ To get an upper bound for $N$
we use the estimate \cite{kuhn} 
\[
\sup_{k\geq 1}k^{1/2}\cdot e_{k}\leq C\cdot n^{1/2}\cdot M_{V^{o}}
\]%
which is valid for any convex symmetric body $V\subset \mathbb{R}^{n}.$ Put $%
e_{k}=1,$ then minimal cardinality $N$ of $1-$net for $J^{-1}\left(
B_{q}\cap \mathcal{T}_{N}\right) $ in the norm induced by $J^{-1}\left(
B_{2}\cap \mathcal{T}_{N}\right) $ can be estimated as 
\[
N\leq 2^{CnM_{V^{o}}^{2}}
\]%
where $V:=J^{-1}\left( B_{q}\cap \mathcal{T}_{N}\right) .$ Finaly, we get%
\[
\mathrm{Vol}_{n}((X_{\mathbf{l}}^{m}+\,J^{-1}(\left( \vartheta _{m}\cdot
B_{q}\right) \cap \mathcal{T}_{2N})\cap J^{-1}(B_{2}\cap \mathcal{T}_{N}))
\]%
\[
\leq 2^{CnM_{V^{o}}^{2}}\cdot \vartheta _{m}^{l}\cdot \mathrm{Vol}%
_{s}\,\left( B_{(2)}^{s}\right) \cdot \mathrm{Vol}_{l}\left(
B_{(2)}^{l}\right) 
\]%
Observe that the number of terms in (\ref{b}) is 
\[
\leq \sum_{m=0}^{\dim \mathcal{T}_{2N}}\left( 
\begin{array}{c}
\dim \mathcal{T}_{2N} \\ 
m%
\end{array}%
\right) =2^{\dim \mathcal{T}_{2N}}.
\]%
Hence, from (\ref{b}) it follows that 
\begin{equation}
\mathrm{Vol}_{n}(J^{-1}(B_{p}\cap \mathcal{T}_{N}))\leq 2^{n}\cdot \max_{%
\mathbf{l}}\,\omega _{m}^{\mathbf{l}},  \label{kukareku}
\end{equation}%
where 
\[
\omega _{m}^{\mathbf{l}}:=\mathrm{Vol}_{n}\left( (X_{\mathbf{l}%
}^{m}+\,J^{-1}(\left( \vartheta _{m}\cdot B_{q}\right) \cap \mathcal{T}%
_{2N})\cap J^{-1}(B_{2}\cap \mathcal{T}_{N})\right) 
\]%
\[
\leq \vartheta _{m}^{l}\cdot 2^{CnM_{\left( J^{-1}\left( B_{q}\cap (JX_{%
\mathbf{l}}^{m})_{\mathcal{T}_{N}}^{\perp }\right) \right) ^{o}}^{2}}
\]%
\begin{equation}
\cdot \mathrm{Vol}_{s}\,\left( B_{(2)}^{s}\right) \cdot \mathrm{Vol}%
_{l}\,\left( J^{-1}\left( B_{(2)}\cap (JX_{\mathbf{l}}^{m})_{\mathcal{T}%
_{N}}^{\perp }\right) \right) .  \label{leic-1}
\end{equation}%
Let $\alpha \in (JX_{\mathbf{l}}^{m})_{\mathcal{T}_{N}}^{\perp }.$ By H\"{o}%
lder's inequality%
\[
\left\Vert \alpha \right\Vert _{\left( q\right) }^{o}=\sup_{\left\Vert \beta
\right\Vert _{\left( q\right) }\leq 1}\left\langle \alpha ,\beta
\right\rangle =\sup_{\left\Vert \beta \right\Vert _{\left( q\right) }\leq
1}\int_{\mathbb{S}^{d}\left( \mathbb{C}\right) }J\alpha \cdot J\beta \cdot
d\nu 
\]%
\[
\leq \left\Vert J\alpha \right\Vert _{q^{^{\prime }}}\cdot \left\Vert J\beta
\right\Vert _{q}
\]%
\[
=\left\Vert J\alpha \right\Vert _{q^{^{\prime }}}=\left\Vert \alpha
\right\Vert _{\left( q^{^{\prime }}\right) },
\]%
Remind that $l=\dim (JX_{\mathbf{l}}^{m})_{\mathcal{T}_{N}}^{\perp }.$ The
last inequality implies%
\[
M_{\left( J^{-1}\left( B_{q}\cap (JX_{\mathbf{l}}^{m})_{\mathcal{T}%
_{N}}^{\perp }\right) \right) ^{o}}=\int_{\mathbb{S}^{l-1}}\left\Vert \alpha
\right\Vert _{\left( q\right) }^{o}\cdot d\mu _{l}
\]%
\[
\leq \int_{\mathbb{S}^{l-1}}\left\Vert \alpha \right\Vert _{\left(
q^{^{\prime }}\right) }\cdot d\mu _{l}=M_{J^{-1}\left( B_{q^{^{\prime
}}}\cap (JX_{\mathbf{l}}^{m})_{\mathcal{T}_{N}}^{\perp }\right) }
\]%
\begin{equation}
\leq C\cdot (q^{^{\prime }})^{1/2},1/q+1/q^{^{\prime }}=1,1<q<\infty ,
\label{21}
\end{equation}%
Where in the last line we use Lemma 3 and (\ref{infty}). Comparing (\ref%
{kukareku})-(\ref{21}) we get 
\begin{equation}
\omega _{m}^{\mathbf{l}}\leq \vartheta _{m}^{^{l}}\cdot 2^{C\cdot
(q^{^{\prime }})^{1/2}\cdot n}\cdot \mathrm{Vol}_{s}\,\left(
B_{(2)}^{s}\right) \cdot \mathrm{Vol}_{l}\,\left( B_{(2)}^{l}\right) .
\label{0}
\end{equation}%
Now we turn to the lower bounds for $\mathrm{Vol}_{n}(J^{-1}(B_{p}\cap 
\mathcal{T}_{N})).$ From the Bourgain-Milman inequality \cite{burm} 
\[
\left( \frac{\mathrm{Vol}_{n}V\cdot \mathrm{Vol}_{n}V^{o}}{\left( \mathrm{Vol%
}_{n}B_{(2)}^{n}\right) ^{2}}\right) ^{1/n}\geq C,
\]%
which is valid for any convex symmetric body $V\subset \mathbb{R}^{n},$ it
follows that 
\[
\mathrm{Vol}_{n}(J^{-1}(B_{p}\cap \mathcal{T}_{N}))\geq C^{n}\cdot \left( 
\frac{\mathrm{Vol}_{n}\left( B_{(2)}^{n}\right) }{\mathrm{Vol}_{n}\left(
B_{(p)}^{n}\right) ^{o}}\right) \cdot \mathrm{Vol}_{n}\left(
B_{(2)}^{n}\right) .
\]%
Comparing this estimate with Lemma 3 and Urysohn's inequality \cite{pisier} 
\[
\left( \frac{\mathrm{Vol}_{n}\left( V\right) }{\mathrm{Vol}_{n}\left(
B_{(2)}^{n}\right) }\right) ^{1/n}\leq \int_{\mathbb{S}^{n-1}}\left\Vert
\alpha \right\Vert ^{o}d\mu ,\,\,\,\Vert \cdot \Vert =\Vert \cdot \Vert _{V},
\]%
which is valid for any convex symmetric body $V\subset \mathbb{R}^{n}$ we
get 
\[
\left( \frac{\mathrm{Vol}_{n}\left( B_{(2)}^{n}\right) }{\mathrm{Vol}%
_{n}\left( B_{(p)}^{n}\right) ^{o}}\right) \cdot \mathrm{Vol}_{n}\left(
B_{(2)}^{n}\right) 
\]%
\[
\geq C^{n}\cdot \left( M_{B_{(p)}^{n}}\right) ^{-n}\mathrm{Vol}_{n}\left(
B_{(2)}^{n}\right) 
\]%
\begin{equation}
\geq C^{n}\cdot \left\{ 
\begin{array}{cc}
(p)^{1/2}, & p<\infty  \\ 
(\log N)^{1/2}, & p=\infty 
\end{array}%
\right\} ^{-n}\cdot \mathrm{Vol}_{n}\left( B_{(2)}^{n}\right) .  \label{zina}
\end{equation}

\bigskip Applying (\ref{kukareku}), (\ref{0}), (\ref{zina}) we obtain%
\[
C^{n}\cdot \left\{ 
\begin{array}{cc}
(p)^{1/2}, & p<\infty  \\ 
(\log N)^{1/2}, & p=\infty 
\end{array}%
\right\} ^{-n}\cdot \mathrm{Vol}_{n}\left( B_{(2)}^{n}\right) 
\]%
\[
\leq 2^{n}\cdot \max_{\mathbf{l}}\vartheta _{m}^{^{l}}\cdot 2^{C\cdot
(q^{^{\prime }})^{1/2}\cdot n}\cdot \mathrm{Vol}_{s}\,\left(
B_{(2)}^{s}\right) \cdot \mathrm{Vol}_{l}\,\left( B_{(2)}^{l}\right) 
\]%
which means that 
\[
2^{-n}\cdot 2^{C\cdot (q^{^{\prime }})^{1/2}\cdot n}\cdot C^{n}\cdot \left\{ 
\begin{array}{cc}
(p)^{1/2}, & p<\infty  \\ 
(\log N)^{1/2}, & p=\infty 
\end{array}%
\right\} ^{-n}
\]%
\[
\times \left( \frac{\mathrm{Vol}_{n}\left( B_{(2)}^{n}\right) }{\mathrm{Vol}%
_{s}\,\left( B_{(2)}^{s}\right) \cdot \mathrm{Vol}_{l}\,\left(
B_{(2)}^{l}\right) }\right) \leq \max_{\mathbf{l}}\vartheta _{m}^{^{l}}.
\]%
or%
\[
2^{-n/l}\cdot 2^{C\cdot (q^{^{\prime }})^{1/2}\cdot n/l}\cdot C^{n/l}\cdot
\left\{ 
\begin{array}{cc}
(p)^{1/2}, & p<\infty  \\ 
(\log N)^{1/2}, & p=\infty 
\end{array}%
\right\} ^{-n/l}
\]%
\begin{equation}
\times \left( \frac{\mathrm{Vol}_{n}\left( B_{(2)}^{n}\right) }{\mathrm{Vol}%
_{s}\,\left( B_{(2)}^{s}\right) \cdot \mathrm{Vol}_{l}\,\left(
B_{(2)}^{l}\right) }\right) ^{1/l}\leq \vartheta _{m}.  \label{zzzz}
\end{equation}%
Observe that this lower bound holds for any $l$. Let, in particular, $%
m:=\dim \mathcal{T}_{[N/2]}$. Observe that $0\leq \dim JX_{\mathbf{l}}\cap
J^{-1}\mathcal{T}_{N}\leq \dim \mathcal{T}_{[N/2]}$ or $\dim \mathcal{T}%
_{[N/2]}\leq (\dim JX_{\mathbf{l}}\cap \mathcal{T}_{N})^{\perp }\leq \dim 
\mathcal{T}_{N}.$ It implies that $\dim \mathcal{T}_{N/2}\leq l\leq \dim 
\mathcal{T}_{N}$ or $Cn\leq l\leq n,$ where $0<C\leq 1.$ Let us put for
convenience $\mathrm{Vol}_{0}\left( B_{(2)}^{0}\right) =1.$ Since%
\[
\mathrm{Vol}_{n}\left( B_{(2)}^{n}\right) =\frac{\pi ^{n/2}}{\Gamma \left(
n/2+1\right) }
\]%
and%
\[
\mathrm{\Gamma }\left( z\right) =z^{z-1/2}\cdot e^{-z}\cdot \left( 2\pi
\right) ^{1/2}\cdot \left( 1+O\left( z^{-1}\right) \right) 
\]%
then 
\[
r_{l,s,n}:=\left( \frac{\mathrm{Vol}_{s}\,\left( B_{(2)}^{s}\right) \cdot 
\mathrm{Vol}_{l}\,\left( B_{(2)}^{l}\right) }{\mathrm{Vol}_{n}\left(
B_{(2)}^{n}\right) }\right) ^{1/l}
\]%
\[
=\left( \frac{\Gamma \left( n/2+1\right) \cdot \pi ^{\left( n-l\right)
/2}\cdot \pi ^{l/2}}{\pi ^{n/2}\cdot \Gamma \left( \left( n-l\right)
/2+1\right) \cdot \Gamma \left( l/2+1\right) }\right) ^{1/l}
\]%
\[
=\left( \frac{\Gamma \left( n/2+1\right) }{\Gamma \left( \left( n-l\right)
/2+1\right) \cdot \Gamma \left( l/2+1\right) }\right) ^{1/l}
\]%
\[
=\left( \frac{e^{-n/2-1}\cdot \left( \frac{n}{2}+1\right) ^{n/2+1-1/2}}{%
e^{-\left( n-l\right) /2-1}\cdot \left( \frac{n-l}{2}+1\right) ^{\left(
n-l\right) /2+1-1/2}\cdot e^{-l/2-1}\cdot \left( \frac{l}{2}+1\right)
^{l/2+1-1/2}}\right) ^{1/l}
\]%
\[
\times \left( \frac{\left( 1+O\left( \frac{1}{n}\right) \right) }{\left(
1+O\left( \frac{1}{n-l}\right) \right) \cdot \left( 1+O\left( \frac{1}{l}%
\right) \right) }\right) ^{1/l}
\]%
\[
\leq C\cdot \left( \frac{\left( \frac{n}{2}+1\right) ^{n/2+1-1/2}}{\left( 
\frac{n-l}{2}+1\right) ^{\left( n-l\right) /2+1-1/2}\cdot \left( \frac{l}{2}%
+1\right) ^{l/2+1-1/2}}\right) ^{1/l}
\]%
\[
=C\cdot \left( \frac{\left( n+2\right) ^{n/2+1-1/2}}{\left( n-l+2\right)
^{\left( n-l\right) /2+1-1/2}\cdot \left( l+2\right) ^{l/2+1-1/2}}\right)
^{1/l}
\]%
\[
\leq C\cdot \left( \frac{n^{n/2+1/2}}{\left( n-l\right) ^{\left( n-l\right)
/2+1/2}\cdot l^{l/2+1/2}}\right) ^{1/l}
\]%
\[
\leq C\cdot \frac{n^{n/\left( 2l\right) +1/\left( 2l\right) }}{\left(
n-l\right) ^{\left( n-l\right) /\left( 2l\right) +1/\left( 2l\right) }\cdot
l^{1/2+1/\left( 2l\right) }}
\]%
\[
\leq C\cdot \frac{n^{n/\left( 2l\right) }}{\left( n-l\right) ^{\left(
n-l\right) /\left( 2l\right) }\cdot l^{1/2}}
\]%
\[
\leq C\cdot \frac{n^{n/\left( 2l\right) }}{\left( n-l\right) ^{n/\left(
2l\right) -1/2}\cdot l^{1/2}},
\]%
where the penultimate and ultimate steps are justified by the condition $%
Cn\leq l<n$ (see, e.g., \cite{kupa}). Consequently, $r_{Cn,n-Cn,n}\leq C$
for any $n\in \mathbb{N}$ and using (\ref{zzzz}) we get 
\[
\vartheta _{m}\geq C\left\{ 
\begin{array}{cc}
p^{1/2}\cdot 2^{C(q^{^{\prime }})^{1/2}}, & p<\infty ,q>1, \\ 
(\log m)^{1/2}\cdot 2^{C(q^{^{\prime }})^{1/2}}, & p=\infty ,q>1,%
\end{array}%
\right\} ^{-C}.
\]%
Finally, from (\ref{00}) and the last line it follows 
\[
\nu _{m}\geq C\cdot m^{-\gamma /d}\cdot \vartheta _{m}.
\]%
$\blacksquare $

\textbf{Remark 1} \emph{Comparing Theorem 1 and Theorem 2 we get 
\[
\nu _{m}(W_{p}^{\gamma },\Xi ,L_{q})\asymp m^{-\gamma /d},\,\,\gamma
>(d-1)/2,\,\,1<q\leq p<\infty . 
\]%
} \textbf{Remark 2} \emph{By embedding we get 
\[
\nu _{m}(W_{p}^{\gamma },\Xi ,L_{q})\gg m^{-\gamma /d},\,\,\gamma
>0,\,\,1<q,p<\infty , 
\]%
\[
\nu _{m}(W_{\infty }^{\gamma },\Xi ,L_{q})\gg m^{-\gamma /d}(\log
m)^{-C},\,\,\gamma >0\,\,1<q<\infty , 
\]%
\[
\nu _{m}(W_{\infty }^{\gamma },\Xi ,L_{1})\gg m^{-\gamma /d-\epsilon }(\log
m)^{-C},\,\,\gamma >0, 
\]%
where $C>0$ is an absolute constant and $\epsilon $ is an arbitrary positive
number. }

\textbf{Acknowledgement} We would like to thank the referees and the
Communicating Editor for the useful suggestions and comments.



\begin{thebibliography}{99}
\bibitem{burm} Bourgain, J., Milman, V. D., New volume ratio properties for
convex symmetric bodies in $\mathbb{R}^{n}$, \emph{Invent. math. }, \textbf{%
88} (1987), 319--340.

\bibitem{dz} Ditzian, Z., A Kolmogorov-type inequality, \emph{Math. Proc.
Camb. Phil. Soc.}, \textbf{136} (2004), 657-663.

\bibitem{gine} Gin\'e, E., The addition formula for the eigenfunctions of
the Laplacian, \emph{Advances in Mathematics} \textbf{18} (1975), 102--107.

\bibitem{ismagilov} Ismagilov, R. S., $n$-Widths of sets in linear normed
spaces and approximation of functions by trigonometric polynomials, \emph{%
Uspekhi Mat. Nauk.} \textbf{29} (1977), 161-178.

\bibitem{kashin} Kashin, B., Tzafriri, L., Lower estimates for the supremum
of some random processes, II, \emph{East J. Approx.} \textbf{1} (1995),
373--377.

\bibitem{kuhn} K\"uhn, T, $\gamma$-radonifying operators and entropy ideals, 
\emph{Math. Nachr.} \textbf{107} (1982), 53--58.

\bibitem{ku1} Kushpel, A. K., On an estimate of Levy means and medians of
some distributions on a sphere, in \emph{Fourier Series and their
Applications}, Inst. of Math., Kiev, 1992, 49--53.

\bibitem{ku2} Kushpel, A. K., Estimates of Bernstein's widths and their
analogs, \emph{Ukrain. Math. Zh.} \textbf{45} (1993), 54--59.

\bibitem{klw} Kushpel, A. K., Levesley, J., Wilderotter, K., On the
asymptotically optimal rate of approximation of multiplier operators from $%
L_{p}$ into $L_{q}$, \emph{J. Constr. Approx.} \textbf{14} (1998), p.
169--185.

\bibitem{ku3} Kushpel, A. K., Levy means associated with two-point
homogeneous spaces and applications, in Annals of the $49^{o}$ Semin\'ario
Brasileiro de An\'alise, 1999, 807--823.

\bibitem{ku5} Kushpel, A. K., $n$-Widths of Sobolev's Classes on Compact
Globally Symmetric Spaces of Rank 1, in Trends in Approximation Theory,
2001, Vanderbilt Univ. Press, Nashville, TN, K. Kopotun, T. Lyche, M. Neamtu
(eds.), 201--210.

\bibitem{kuto1} Kushpel, A. K., Tozoni, S. A., On the Problem of Optimal
Reconstruction, \emph{Journal of Fourier Analysis and Applications} \textbf{%
13} (2007), 459--475.

\bibitem{klt2009} Kushpel, A. K., Levesley, J., Tozoni, S. A., Estimates of $%
n$-widths of Besov classes on two-point homogeneous manifolds, \emph{%
Matematische Nachrichten} \textbf{282} (2009), 748--763. 1621--1629.

\bibitem{kupa} Kushpel, A., Optimal cubature formulas on compact homogeneous
manifolds, \emph{Journal of Functional Analysis} \textbf{257} (2009),
1621--1629.

\bibitem{lk} Kushpel, A., Levesley, J., A Multiplier version of the
Bernstein inequality on the complex sphere, \emph{Technical Report,
University of Leicester} \textbf{257} (2011), 1--16.

\bibitem{kwap} Kwapie\'n, S., Isomorphic characterizations of inner product
spaces by orthogonal series with vector valued coefficients, \emph{Studia
Math.} \textbf{44}, 583--595.

\bibitem{pisier} Pisier, G., \emph{The volume of convex bodies and Banach
space geometry}, Cambridge Univ. Press, London, 1989.

\bibitem{stechkin} Stechkin, S., On absolute convergence of orthogonal
series, \emph{Dokl. Akad. Nauk SSSR (N.S.)}, 102 (1955) 37--40. (in Russian).

\end{thebibliography}
\end{document}